\documentstyle[amssymb]{article}

\font\boldscriptfont=cmbx7

\font\boldscriptscriptfont=cmbx5

\let\oldbf=\bf
\def\newbf#1{\ifmmode\mathchoice
{\hbox{\oldbf #1}}
{\hbox{\oldbf #1}}
{\hbox{\boldscriptfont #1}}
{\hbox{\boldscriptscriptfont #1}}
\else\oldbf #1\fi}

\def\Bbb{\ifmmode\let\next\Bbb@\else
 \def\next{\errmessage{Use \string\Bbb\space only in math mode}}\fi\next}
\def\Bbb@#1{{\Bbb@@{#1}}}
\def\Bbb@@#1{\fam\msyfam#1}

\reversemarginpar

\reversemarginpar

\newcommand{\lbv}{[\![}
\newcommand{\rbv}{]\!]}

\newcommand{\rest}{{\mathord{\restriction}}}
\newcommand{\add}{{\mathrm {\bold {add}}}\/}
\newcommand{\cov}{{\mathrm {\bold {cov}}}\/}
\newcommand{\unif}{{\mathrm {\bold {unif}}}\/}
\newcommand{\cof}{{\mathrm {\bold {cof}}}\/}
\newcommand{\cf}{\mathrm{cf}}
\newcommand{\dom}{{\mathrm {dom}}}

\newcommand{\QED}{\vrule width 6pt height 6pt depth 0pt \vspace{0.1in}}
\newcommand{\forces}{\mathrel{\|}\joinrel\mathrel{-}}

\newcommand{\N}{{\cal N}}
\newcommand{\M}{{\cal M}}
\newcommand{\ran}{{\mathrm { ran}}}
\newcommand{\V}{{\bold V}}
\newcommand{\<}{\langle}
\renewcommand{\>}{\rangle}

\newcommand{\thinks}{\models}
\newtheorem{theorem}{Theorem}[section]
\newtheorem{lemma}[theorem]{Lemma}

\newtheorem{definition}[theorem]{Definition}
\newcommand{\lesdot}{\mathrel{\mathord{<}\!\!\raise 
0.8 pt\hbox{$\scriptstyle\circ$}}}
\newcommand{\Proof}{{\sc Proof} \hspace{0.2in}}

\newcommand{\lft}[2]{\mathopen\ifcase#1{}\oo\or
                        \big#2\or\Big#2\else\oo\fi} 
\newcommand{\rgt}[2]{\mathclose\ifcase#1{}\oo\or
                        \big#2\or\Big#2\else\oo\fi}

\begin{document}

\title{Closed measure zero sets}
\author{Tomek Bartoszynski\thanks{The author thanks the  Lady Davis
Fellowship 
Trust for full support} \\ 
Boise State University\\
Boise, Idaho\\
and\\
Hebrew University \\
Jerusalem
\and
Saharon Shelah\thanks{Research partially supported by Basic Research
Fund, Israel Academy of Sciences, publication 439} \\
Hebrew University\\
Jerusalem}
\maketitle
\begin{abstract}
We study the relationship between the $\sigma$-ideal generated by
closed measure zero sets and the ideals of null and meager sets. 
We show that the additivity of the ideal of closed measure zero sets
is not bigger than covering for category. As a consequence we get
that the additivity of the ideal of closed measure zero sets is
equal to the additivity of the ideal of meager sets.
\end{abstract}
\section{Introduction}
Let $\M$ and $\N$ denote the ideals of meager and null 
subsets of $2^\omega$ respectively and let $
{\cal E}$ be the $\sigma$-ideal generated by closed measure zero
subsets of $2^\omega$. It is clear that ${\cal E}$ is a proper subideal
of $\M \cap \N$.

For an ideal ${\cal J}$  of subsets of $2^\omega$ define
\begin{enumerate}
\item $\add({\cal J})=\min\{|{\cal A}| : {\cal A} \subseteq {\cal J} \ \& \
\bigcup {\cal A} \not \in {\cal J} \}$,
\item $\cov({\cal J}) = \min\{|{\cal A}| : {\cal A} \subseteq {\cal J} \ \&
\ \bigcup {\cal A}= 2^\omega\}$,
\item $\unif({\cal J}) = \min\{|X| : X  \subseteq 2^\omega  \ \&
\ X \not \in {\cal J}\}$ and
\item $\cof({\cal J})=\min\{|{\cal A}| : {\cal A} \subseteq {\cal J} \ \& \
\forall B \in {\cal J} \ \exists A \in {\cal A} \ B \subseteq A \}$.
\end{enumerate}
We can further generalize these definitions and put
for a pair of ideals $ {\cal I} \subseteq {\cal J}$,
\begin{enumerate}
\item $\add({\cal I},{\cal J})=\min\{|{\cal A}| : {\cal A} \subseteq {\cal
I} \ \& \ \bigcup {\cal A} \not \in {\cal J}\}$,
\item $\cof({\cal I},{\cal J})=\min\{|{\cal A}| : {\cal A} \subseteq {\cal
J} \ \& \ \forall B \in {\cal I} \ \exists A \in {\cal A} \ B
\subseteq A\}$.
\end{enumerate}

Let ${\cal I}_0$ be the ideal of finite subsets of $2^\omega$. Note
that $\cov({\cal J})=\cof({\cal I}_0, {\cal J})$, $\unif({\cal
J})=\add({\cal I}_0, {\cal J})$, $\add({\cal J})=\add({\cal
J},{\cal J})$ and $\cof({\cal J})=\cof({\cal J},{\cal J})$.

The goal of this paper is to study the relationship between the
cardinals defined above for the ideals $\M$,  $\N$ and ${\cal E}$. We
will show that $\add(\M)=\add({\cal E})$ and $\cof(\M)=\cof({\cal
E})$.

It will follow from the inequalities
$\add({\cal E},\N) \leq  \cov(\M)$ and $\cof( {\cal E}, \N) \geq
\unif(\M)$ which will be proved in section 3.

Finally in the last section we will present some consistency results --
we will show the $\cov({\cal E})$ may not be equal to $\max\{\cov(\N),
\cov(\M)\}$ and similarly $\unif({\cal E})$ does not have to be equal
to $\min\{\unif(\M), \unif(\N)\}$.
 
For $f,g \in \omega^\omega$ let $f \leq^\star g$ be the ordering of
eventual dominance.

Recall that ${\frak b}$ is the size of
the smallest unbounded family in $\omega^\omega$ and ${\frak d}$ is the
size of the smallest dominating family in $\omega^\omega$.

Through this paper we use the standard notation. 

$\mu$ denotes the
standard product measure on $2^\omega$. 
For a  tree $T \subseteq 2^{<\omega}$ let
$[T]$ be the set of branches of $T$. If $T$ is finite (or has terminal
nodes) then $[T]$
denotes the clopen subset of $2^\omega$ determined by maximal nodes of
$T$. Let $m(T)=\mu([T])$ in both cases. 

If $s \in T \subseteq
2^{<\omega}$ then $T[s]=\{t : s^\frown t \in T\}$ where $s^\frown t$
denotes the concatenation of $s$ and $t$.
ZFC$^\star$ always denotes some finite fragmet of ZFC sufficiently big
for our purpose.

We will conclude this section with several results concerning the
cardinal invariants defined above.

\begin{theorem}[Miller {[Mi]}]\label{miller}
\ 
\begin{enumerate}
\item $\add(\M)=\min\{\cov(\M), {\frak b}\}$ and
$\cof(\M)=\max\{\unif(\M), {\frak d}\}$,
\item $\add({\cal E},\M) \leq {\frak b}$ and $\cof({\cal E},\M) \geq
{\frak d}$. In particular $\add({\cal E}) \leq {\frak b}$ and
$\cof({\cal E}) \geq {\frak d}$,
\item $\cov(\M) \leq \add({\cal E},\N)$ and $\unif(\M)
\geq \cof({\cal E},\N)$.~$\QED$
\end{enumerate}
\end{theorem}

We will also use the combinatorial characterizations of cardinals
$\cov(\M)$ and $\unif(\M)$.

\begin{theorem}[Bartoszynski {[Ba1]}]\label{bart}
\ 
\begin{enumerate}
\item $\cov(\M)$ is the size of the smallest family $F \subseteq
\omega^\omega$ such that 
$$\forall g \in \omega^\omega \ \exists f \in F \ \forall^\infty n \
f(n) \neq g(n).$$
\item $\unif(\M)$ is the size of the smallest family $F \subseteq
\omega^\omega$ such that 
$$\forall g \in \omega^\omega \ \exists f \in F \ \exists^\infty n \
f(n)=g(n). \ \QED$$
\end{enumerate}
\end{theorem}

\section{Combinatorics}
\setcounter{theorem}{0}
In this section we will prove several combinatorial lemmas which will
be needed later. The following theorem uses the technique from [Ba2].
\begin{theorem}\label{smalunion}
Suppose that $\{F_\eta : \eta < {\lambda} <\add({\cal E},\N)\}$ is a
family of closed measure zero sets.
Then there exists a partition of $\omega$ into intervals $\{\bar{I}_n : n
\in \omega\}$ and a sequence $\{T_n : n \in \omega\}$ such that 
for all $n$, $T_n \subseteq 2^{\bar{I}_n}$, $|T_n| \cdot 2^{-|\bar{I}_n|} \leq
2^{-n}$ and
$$\bigcup_{\eta < {\lambda}} F_\eta \subseteq \{x \in 2^\omega :
\exists^\infty n \ x \rest \bar{I}_n  \in T_n\}.$$
Furthermore, we can require that
$$\forall \eta<{\lambda} \ \exists^\infty n \ F_\eta \rest \bar{I}_n
\subseteq T_n$$
where $F_\eta \rest \bar{I}_n = \{s \in 2^{\bar{I}_n} : \exists x \in
F_\eta \ x \rest \bar{I}_n = s\}$.
\end{theorem}
\Proof
Note that if the sequences $\{\bar{I}_n : n \in \omega\}$ and $\{T_n : n \in
\omega\}$ satisfy the above conditions then the set 
$\{x \in 2^\omega : \exists^\infty n \ x \rest \bar{I}_n \in T_n\}$ has
measure zero.

For $\eta<{\lambda}$ and $n \in \omega$ define
$$F^n_\eta = \{x \in 2^\omega : \exists s \in 2^n \ s^\frown x \rest
(\omega-n) \in F_\eta\}.$$

By the assumption there exists a measure zero set $H \subseteq
2^\omega$ such that  $\bigcup_{\eta < {\lambda}} \bigcup_{n \in \omega}
F^n_\eta \subseteq H$.

\begin{lemma}[Oxtoby {[O]}]\label{oxtoby} 
There exists a sequence of finite sets $\< H_n : n  \in
\omega  \>$  such
that  $H_n  \subseteq  2^{n}$,
$\sum_{n=1}^{\infty}   |H_n|\cdot2^{-n} <
\infty$   and  $H  \subseteq  \{x  \in  2^{\omega}  :
\exists^{\infty} n \ x
\rest n  \in  H_n \}$.
\end{lemma}
\Proof Since $H$ has measure zero there are open sets
$\<G_{n} : n \in \omega \>$ 
covering $H$ such that $\mu(G_{n}) < 2^{-n}$
for $n \in \omega$.
Represent each set $G_{n}$ as a disjoint union of open basic intervals 
$$G_{n} = \bigcup_{m=1}^{\infty} [s^{n}_{m}]
\hbox{ for } n \in \omega  .$$
Let $H_n = \{ s \in 2^{n} : s = s^{k}_{l} $ for some $k,l \in \omega \}$
for $n \in \omega$. It follows that $\sum_{n=1}^\infty |H_n| \cdot
2^{-n} \leq \sum_{n=1}^\infty \mu(G_n) \leq 1$.
If $x \in H$ then $x \in \bigcap_{n \in \omega} G_n$. Therefore $x
\rest n \in F_n$ must hold for infinitely many $n$.~$\QED$

Therefore 
$$\bigcup_{\eta<{\lambda}} \bigcup_{n \in \omega} F^n_\eta 
\subseteq \{x \in 2^\omega  :
\exists^\infty n \ x \rest n \in H_n\}.$$
For every $\eta<{\lambda}$ define an increasing sequence $\<k^\eta_n :
n \in \omega\>$ as follows: $k^\eta_0=0$ and for $n \in \omega$,
$$k^\eta_{n+1} = \min\left\{m : F^{k^\eta_n}_\eta \subseteq
\bigcup_{j=k^\eta_n}^m [H_j] \right\}.$$
Since sets $F^n_\eta$ are compact this definition is correct.

We will need an increasing sequence $\<k_n : n \in \omega\>$ such that
$$\forall \eta<{\lambda} \ \exists^\infty n \ \exists m \ k_{2n}
<k^\eta_m <k^\eta_{m+1} < k_{2n+1}$$
and
$$2^{k_n} \cdot \sum_{j=k_{n+1}}^\infty \frac{|H_j|}{2^j} \leq
\frac{1}{2^n}.$$
To construct such a sequence we will use the following lemma:
\begin{lemma}
Suppose that $M \thinks {\mathrm {ZFC}}^\star$ and $|M|<{\frak d}$.
Then there exists a function $g \in \omega^\omega$ such that either
$$\forall f \in M \cap \omega^\omega \ \exists^\infty n \ \exists m \
g(2n)<f(m)< f(m+1)<g(2n+1)$$
or
$$\forall f \in M \cap \omega^\omega \ \exists^\infty n \ \exists m \
g(2n+1)<f(m)< f(m+1)<g(2n+2).$$
\end{lemma}
\Proof
Let $g \in \omega^\omega$ be an increasing function such that $g \not
\leq^\star f$ for $f \in M \cap \omega^\omega$.
We will show that $g$ has required properties.

Suppose not. Let $f_1, f_2 \in M \cap \omega^\omega$ be such that 
for all $n$,
$$|[g(2n),g(2n+1)] \cap \ran(f_1)|\leq 1 \hbox{ and } |[g(2n+1),g(2n+2)] \cap
\ran(f_2)|\leq 1.$$
We will get a contradiction by constructing 
a function $f \in M \cap \omega^\omega$ which
dominates $g$.

Define $f(0)=f_1(0)>g(0)$ and $f(1)=f_2(0)>g(1)$.
Let $l_1= \min\{l : f_1(l)>f_2(1)\}$ and put
$f(2)=f_1(l_1)$. Now $f(2)>g(2)$ since $f_2(1)>g(2)$.
Let $l_2=\min\{l : f_2(l)>f_1(l_1+1)\}$ and let $f(3)=f_2(l_2)>g(3)$
since $f_1(l_1+1)>g(3)$. And so on \dots.

In general define the sequence $\<l_n : n \in \omega\>$ as $l_0=0$ and
$$l_{2n+1}=\min\{l : f_1(l)>f_2(l_{2n}+1)\}$$
and
$$l_{2n+2}=\min\{l : f_2(l) > f_1(l_{2n+1}+1)\}.$$
Let 
$$f(n+1)= \left\{\begin{array}{ll}
f_1(l_n) & \hbox{if } n \hbox{ is even }\\
f_2(l_n) & \hbox{if } n \hbox{ is odd}
\end{array}\right. . $$
It is clear that $f \in M$. Easy induction shows that $f $ dominates
$g$. Contradiction.~$\QED$

To get the sequence the desired sequence 
$\<k_n : n \in \omega\>$ take a
model $M \thinks {\mathrm {ZFC}}^\star$ containing $\<H_n : n \in \omega\>$
and $\{F_\eta : \eta < {\lambda}\}$. Since ${\lambda} < \add({\cal
E},\N) \leq {\frak d}$ we can assume that $|M|<{\frak d}$.
Apply the above lemma to get a function $g$ and define $k_n=g(n)$ for
$n \in \omega$. It is clear that this is the sequence we are looking
for.

Now define for $n \in \omega$,
$$\bar{I}_n = [k_{2n-1},k_{2n+1}]$$ and
$$T_n = \left\{s \in 2^{\bar{I}_n} : \exists j \in [k_{2n},k_{2n+1}] \
\exists t \in H_j \ s \rest \bar{I}_n = t \rest \bar{I}_n \right\}.$$
Note that for every $n$,
$$\frac{|T_n|}{2^{|\bar{I}_n|}} \leq 2^{k_n} \cdot \sum_{j=k_{2n}}^{k_{2n+1}}
\frac{|H_j|}{2^j} \leq \frac{1}{2^n} .$$

To finish the proof fix $\eta < {\lambda}$ and $k \in \omega$.
By the construction there exists $n>k$ and $m \in \omega$ such that 
$$k_{2n}<k^\eta_m<k^\eta_{m+1}<k_{2n+1}.$$
Suppose that $s \in F_\eta \rest \bar{I}_n$.
Then there exists $x \in F^{k^\eta_m}_\eta$ such that $s \subseteq x$.
Furthermore, there exists $j \in [k^\eta_m, k^\eta_{m+1})$ such that
$x \rest j \in H_j$. It follows that $s \in T_n$.~$\QED$

Now we will prove another combinatorial lemma describing the structure
of closed measure zero sets.

Let $\left\{I_n :n \in \omega\right\}$ be a partition of $\omega$ into
disjoint intervals such that $|I_n|>n$.

For $n<m$ let 
$$Seq_{n,m}=\left\{s : \dom(s) \subseteq [n,m] \ \& \ \forall j \in \dom(s)
\ s(j) \in I_j\right\}.$$
For every $s \in Seq_{n,m}$ define
$$C_s = \left\{t : \dom(t)=\bigcup_{j=n}^m I_j \ \& \ \forall j \in
\dom(s) \ 
t \lft1( s(j)\rgt1) =0\right\}.$$ 
For $k,j \in \omega $ let 
$$C^j_{k}=\left\{\begin{array}{ll}
\left\{t \in 2^{I_j} : t(k)=0\right\} & \hbox{if } k \in I_j\\
2^{I_j} & \hbox{otherwise}
\end{array}\right. .$$
Note that we can identify the set $C_s$ with $\prod_{j=n}^m
C^j_{s(j)}$ in the following way:
$$t \in C_s \ \leftrightarrow \ \exists \<t_n,t_{n+1}, \ldots, t_m\>
\in
\prod_{j=n}^m C^j_{s(j)} \ \ t = t_n\!^\frown t_{n+1}\!^\frown
\cdots^\frown t_m .$$

Fix $n<m$ and let $I =  I_n \cup I_{n+1} \cup \cdots \cup I_m$. 
Suppose that $T \subseteq
2^I$ is a finite tree such that 
\begin{enumerate}
\item $\forall s \in T \ \exists t \in T \ (s \subseteq t \ \& \
|t|=|I|)$,
\item $m(T) \leq \frac{1}{4}$.
\end{enumerate}

\begin{lemma}\label{first}
Suppose that  for some $s \in Seq_{n,m}$, $C_s = \prod_{j=n}^m C^j_{s(j)} 
\subseteq T$. Then there
exists $k \in [n,m)$ and $t \in T \cap \prod_{j=1}^{k-1} C^j_{s(j)}$
(if $k=n$ then $t = \emptyset$) such that 
$$\forall t' \in C^{k}_{s(k)} \ m \lft1( T[t^\frown t']\rgt1)  >
\left(1+\frac{1}{2^k}\right) \cdot m \lft1( T[t] \rgt1) .$$
\end{lemma}
\Proof
Suppose not. We build by induction a sequence $\<t_j : j \in [n,m-1]\>$
such that $t_j \in C^j_{s(j)}$ and 
$m \lft1( T[t_j \!^\frown t_{j+1}] \rgt1)  \leq
(1+2^{-j})\cdot m \lft1( T[t_j] \rgt1) $ for $j<m$.

After $m-1$ many steps we get that
$$m \lft1( T[t_n \!^\frown t_{n+1}\!^\frown \cdots^\frown
t_{m-1}]\rgt1)  \leq m(T)
\cdot \prod_{j=n}^{m-1} \left(1+ \frac{1}{2^j}\right) < \frac{1}{2}.$$
Therefore there is $t_m \in C^m_{s(m)}-
T[t_n \!^\frown t_{n+1}\!^\frown \cdots^\frown t_{m-1}].$
This is a contradiction since
$$t=t_n\!^\frown t_{n+1}\!^\frown \cdots^\frown t_m \in C_s -T. \
\QED$$
Suppose that $t \in T$ and $|t|=|\bigcup_{j=n}^k I_j|$ for some $k \in [n,m)$.
Let
$$S^{k+1}_t = \left\{l \in I_{k+1} : \forall t' \in C^{k+1}_l \
m \lft1( T[t^\frown t'] \rgt1)  > \left(1+\frac{1}{2^k}\right) \cdot m
\lft1( T[t] \rgt1) \right\}.$$
Note that the sets $\left\{C^{k+1}_l : l \in I_{k+1}\right\}$ are independent.
Therefore the set
$$\bigcup_{l \in S^{k+1}_t} \bigcup_{t' \in C^{k+1}_l} T[t^\frown
t']$$
has measure at least
$$\left(1-2^{-|S^{k+1}_t|}\right) \cdot \left(1+\frac{1}{2^k}\right)
\cdot m \lft1( T[t] \rgt1) .$$
Since this set is included in $T[t]$ we get
$$\left(1-2^{-|S^{k+1}_t|}\right) \cdot \left(1+\frac{1}{2^k}\right) \leq 1.$$
Therefore
$$|S^{k+1}_t| \leq k+1.$$

Let $S^{k+1}= \left\{l \in I_{k+1} : 
\exists t \in T \ l \in S^{k+1}_t\right\}$.
Then
$$|S^{k+1}| \leq (k+1) \cdot 
\prod_{j=n}^k 2^{|I_j|} .$$
Also if $t=\emptyset$ then define
$$S^n_\emptyset = \left\{l \in I_n : \forall t' \in C^n_l \ m \lft1(
T[t'] \rgt1) > \left(1+\frac{1}{2^n}\right) \cdot m(T) \right\}.$$
Similarly we get $|S^n_\emptyset| \leq n+1$.

Note that in particular we get that the size of $S^{k}$ does not depend on
the size of $I_{k}$.

Combining \ref{first} with the observations above we get the
following:
\begin{lemma}\label{second}
Suppose that $I=I_n \cup I_{n+1} \cup \cdots \cup I_m$ and $T
\subseteq 2^I$ such that $m(T) < \frac{1}{4}$. Then
there exists a sequence $\<S^k : k \in [n,m]\>$ such that 
\begin{enumerate}
\item $S^k \subseteq I_k $ for $k \in [n,m]$,
\item $|S^k| \leq (k+1) \cdot \prod_{j=n}^{k-1} 2^{|I_j|}$ for $k \in
(n,m]$ and $|S^n| \leq n+1$,
\item for every $s \in Seq_{n,m}$, if $C_s \subseteq T$ then
there exists $k \in [n,m]$ such that $s(k) \in S^k$.~$\QED$
\end{enumerate}
\end{lemma}

We conclude this section with a theorem of Miller which gives an upper
bound for $\cov({\cal E}, \N)$. We will prove it here for
completeness.
\begin{theorem}[Miller {[Mi]}]\label{miller1}
$\add({\cal E},\N) \leq {\frak d}$ and $\cof({\cal E},\N) \geq {\frak
b}$.
\end{theorem}
\Proof
Suppose that $H \subseteq 2^\omega$ is a measure zero set. Using
\ref{oxtoby}, we can find a sequence $\<H_n : n \in \omega\>$ such
that $H_n  \subseteq 2^n$, $\sum_{n=1}^\infty |H_n|\cdot 2^{-n} \leq
\frac{1}{4}$ and
$$H \subseteq \{x \in 2^\omega : \exists^\infty n \ x \rest n \in H_n\}.$$
Define for $n \in \omega$,
$$f_H(n) = \min\left\{m : \sum_{j=m}^\infty \frac{|H_j|}{2^j} <
\frac{1}{4^n}\right\}.$$ 
Suppose that $f \in \omega^\omega$ is an increasing function.
Let
$$G_f = \{x \in 2^\omega : \forall n \ x \lft1( f(n) \rgt1) =0\}.$$
Clearly $G_f$ is a closed measure zero set.
\begin{lemma}
If $f_H \leq^\star f$ then $G_f \not \subseteq H$.
\end{lemma}
\Proof
Suppose that $f_H \leq^\star f$. Without loss of generality we can
assume that $f_H(n) < f(n)$ for all $n$.
For $n \in \omega$ define
$$\bar{H}_n = \left\{ s \in 2^{f_H(n+1)} : \exists j \in
\lft1[f_H(n),f_H(n+1) \rgt1) \ \exists t \in H_j \ s \rest j = t
\right\}.$$
Note that for all $n$,
$$[\bar{H}_n]= \bigcup_{j = f_H(n)}^{f_H(n+1)} [H_j] \hbox{ and }
m(\bar{H}_n) \leq 4^{-n}.$$

By compactness, if $G_f \subseteq H$ then for some $n$,
$$G_f \subseteq \bigcup^{f_F(n+1)}_{j=1} [H_j] = \bigcup_{j \leq n}
[\bar{H}_j] .$$
We will show that this inclusion fails for every $n$ which will give a
contradiction.

Fix $n \in \omega$. Note that it is enough to find $s \in 2^{f_H(n+1)}$
such that $s \lft1( f(j) \rgt1) =0 $  and $s \rest f_H(j+1)
\not \in \bar{H}_j$ for $j \leq n$.

We will use the following simple construction.
\begin{lemma}\label{tricky}
Suppose that $n_1<n_2<n_3$ and  that $T \subseteq 2^{[n_1,n_3]}$ is
such that 
$m(T)~=~a~<~\frac{1}{2}$.
For $l \in [n_2,n_3]$ let $C_l = \left\{s
\in 2^{[n_2,n_3]} : s(l)=0 \right\}$.
Then for every $l \in [n_2,n_3]$ there exists $s \in C_l$ such that 
the set $T[s] = \left\{t \in 2^{[n_1,n_2)} : t^\frown s \in T \right\}$
has measure $\leq 2a$.
\end{lemma}
\Proof
Fix $l \in [n_2,m_3]$ and choose $s$ such that $m(T[s])$ is minimal.

If $T[s] = \emptyset$ we are done.
Otherwise
$$m(T) \geq \frac{1}{2} \cdot m(T[s])  .$$
It follows that $m(T[s]) \leq 2a$.~$\QED$

We will build by induction sequences $s_n, s_{n-1}, \ldots, s_0$ and
sets $H'_n, H'_{n-1}, \ldots, H'_0$ such that for all $j \leq n$,
\begin{enumerate}
\item $\dom(s_j)=[f_H(j),f_H(j+1))$,
\item $H'_j \subseteq 2^{f_H(j+1)}$,
\item $m(H'_j[s_j]) \leq 2\cdot m(H'_j)$.
\end{enumerate}
Let $H'_n=\bar{H}_n$ and let $s_n \in 2^{[f_H(n),f_H(n+1))}$ be the sequence
obtained by applying \ref{tricky} to $H'_n$ and $C_{f(n)}$.

Suppose that $H'_{n-j}$ and $s_{n-j}$ are already constructed.
Let 
$$H'_{n-j-1} = \bar{H}_{n-j-1} \cup H'_{n-j}[s_{n-j}]$$
and let $s_{n-j-1}$ be the sequence obtained by applying
\ref{tricky} to $H'_{n-j-1}$ and $C_{f(n-j-1)}$.

Let $s = s_0\!^\frown s_1\!^\frown \cdots^\frown s_n$.
Note that $s \lft1( f(j)\rgt1) =0$ for all $j \leq n$.
We have to check that $s \rest f_H(j+1) \not \in \bar{H}_j$ for $j
\leq n$.
Suppose this is not true. Pick minimal $j$ such that 
$$s \rest f_H(j+1) = s_0\!^\frown s_1\!^\frown \cdots^\frown s_{j} \in
\bar{H}_{j}.$$  
By the choice of $s_j$ we have
$$s_0\!^\frown s_1\!^\frown \cdots^\frown s_{j-1} \in \bar{H}_{j-1}
\cup \bar{H}_j[s_j].$$
Since $j$ was minimal,
$$s_0\!^\frown s_1\!^\frown \cdots^\frown s_{j-1} \in
\bar{H}_j[s_j].$$
Proceding like that we  get that
$$s_0\!^\frown s_1\!^\frown \cdots^\frown s_{j-2} \in
\bar{H}_j[s_j][s_{j-1}]$$ 
Finally
$$s_0 \in \bar{H}_j[s_j][s_{j-1}]\cdots[s_1] \subseteq H'_0$$
which is a contradiction.~$\QED$

Now we are ready to finish the proof of the theorem.
Suppose that $F \subseteq \omega^\omega$ is a dominating family which
consists of increasing functions.
Consider the set $\bigcup_{f \in F} G_f$. We claim that this set does
not have measure zero. It follows from the fact that if $H$ is a
measure zero set then there exists $f \in F$ such that $f_H
\leq^\star f$.
In particular $G_f \not \subseteq H$.

Similarly, if ${\cal B} \subseteq \N$ is a family of size $< {\frak
b}$ then there exists $f \in \omega^\omega$ such that 
$$\forall H \in {\cal B} \ f_H \leq^\star f .$$
Thus $G_f \not \subseteq H$ for any $H \in {\cal B}$.~$\QED$

\section{Cohen reals from closed measure zero sets}
\setcounter{theorem}{0}
The goal of this section is to prove that $\add({\cal
E},\N)=\cov(\M)$.
In fact we have the following:
\begin{theorem}
\ 
\begin{enumerate}
\item $\add({\cal E},\N)=\cov(\M)$. In particular $\add({\cal
E})=\add(\M)$,
\item $\cof({\cal E},\N)=\unif(\M)$. In particular $\cof({\cal E})=\cof(\M)$.
\end{enumerate}
\end{theorem}
\Proof
Note that by \ref{miller} and \ref{miller1}, we get
$$\add(\M) = \min\{\cov(\M), {\frak b}\} \leq \add({\cal E},\N) \leq
{\frak b}.$$
Therefore the equality $\add({\cal E})=\add(\M)$ follows from the
inequality $\add({\cal E},\N) \leq \cov(\M)$.

Similarly, to show that $\cof({\cal E})=\cof(\M)$ we have to check
that $\cof({\cal E},\N) \geq \unif(\M)$.

\vspace{0.1in}

$(1)$ $\add({\cal E},\N) \leq \cov(\M)$. 

By the first part of \ref{bart}, it is enough to prove that for
every family $F \subseteq \omega^\omega$ of size $<\add({\cal E},\N)$
there exists a function $g \in \omega^\omega$ such that 
$$\forall f \in F \ \exists^\infty n \ f(n)=g(n) .$$
Fix a family $F$ as above.

For every $f \in F$ let
$$f'(n) = \max\left\{f(i) : i \leq n\right\}+1 \hbox{ for } n \in \omega.$$
We will need two increasing sequences $\left\{m_n,l_n : n \in \omega\right\}$
such that 
\begin{enumerate}
\item $m_0=l_0=0$,
\item $l_{n+1} = l_n + 2^{m_n} \cdot
(n+1)$,
\item $\forall f \in F \ \exists^\infty n \ m_{n+1} >
f'(l_{n+1})^{l_{n+1}} + m_n$.
\end{enumerate}
The existence of these sequences follows from the fact that $|F| <
{\frak d}$.

Let $I_n = [m_n,m_{n+1})$ and $J_n = [l_n, l_{n+1})$ for $n \in
\omega$.
Without loss of generality we can assume that $|I_n|={K_n}^{|J_n|}$
for some $K_n \in \omega$.
Thus we can identify elements of $I_n$ with ${K_n}^{J_n}$.

For every $f \in F$ and $n \in \omega$ define $\vec{f}(n)=f \rest
J_n$.
By the choice of sequences $\<I_n, J_n : n \in \omega\>$ we have
$$\forall f \in F \ \exists^\infty n \ \vec{f}(n) \in I_n.$$
Using the notation from previous section, define for $f \in F$,
$$C_f = \bigcap_{n \in \omega} C_{\vec{f} \rest n} .$$
Note that the sets $C_f$ are closed sets of measure zero.

Since $|F|<\add({\cal E},\N)$, the set
$\bigcup_{f \in F} C_f \hbox{ has measure zero.}$

By \ref{smalunion}, there exist sequences $\<\bar{I}_n, T_n : n \in
\omega\>$ such that 
for all $n$, $T_n \subseteq 2^{\bar{I}_n}$, $|T_n|\cdot 2^{-|\bar{I}_n|} \leq
2^{-n}$ and
$$\forall f \in F \ \exists^\infty n \ C_f \rest \bar{I}_n \subseteq
T_n .$$
Moreover, without loss of generality we can assume
that whenever $I_m \cap \bar{I}_n \neq \emptyset$ then $I_m \subseteq
\bar{I}_n$ for $n,m \in \omega$.

We will build the function $g \in \omega^\omega$ 
we are looking for from the sequences
$\<T_n : n~\in~\omega\>$ and $\<I_n : n \in \omega\>$.

For every $n$ let $v_n \in \omega $ be such that 
$$\bar{I}_n = I_{v_n} \cup I_{v_n+1} \cup \cdots \cup I_{v_{n+1}-1} .$$
Note that for $f \in F$ and $n \in \omega$,
$$C_f \rest \bar{I}_n = C_{\vec{f} \rest [v_n, v_{n+1})} .$$
Now we are ready to define function $g$.
For every $n$ we will define $g \rest \bar{I}_n$ using the set $T_n$.

Fix $n \in \omega$ and consider the set $T_n \subseteq 2^{\bar{I}_n}$.
By \ref{second} there exists
a sequence $\<S^k : k \in [v_n, v_{n+1})\>$ such that 
\begin{enumerate}
\item $S^k \subseteq I_k $ for $k \in [v_n,v_{n+1})$,
\item $|S^k| \leq (k+1) \cdot \prod_{j=n}^{k-1} 2^{|I_j|}$ for $k \in
(v_n,v_{n+1})$ and $|S^{v_n}| \leq n+1$,
\item for every $s \in Seq_{v_n,v_{n+1}-1}$, if $C_s \subseteq T$ then
there exists $k \in [v_n,v_{n+1})$ such that $s(k) \in S^k$.
\end{enumerate}

Note that for every $k \in [v_n,v_{n+1})$,
$$|S^k| \leq  (k+1) \cdot \prod_{j=n}^{k-1} 2^{|I_j|} \leq 
(k+1) \cdot \prod_{j=n}^{k-1} 2^{m_{j+1}-m_j} \leq (k+1) \cdot
2^{m_j} \leq |J_k|.$$
 
We can view $S^k$ as a subset of $K_k^{J_k}$ of size $\leq |J_k|$.
For $k \in [v_n,v_{n+1})$ let $s^k \in K_k^{J_k}$ be such that 
$$\forall t \in S^k \ \exists l \in J_k \ s^k(l)=t(l).$$
Define 
$$g \rest \bar{I}_n = s^{v_n}\!^\frown \cdots^\frown s^{v_{n+1}-1}
.$$
Note that $g \rest \bar{I}_n$ ``diagonalizes'' all sets $S^k$ for $k
\in [v_n, v_{n+1})$.

Now we are ready to finish the proof. Suppose that $f \in F$.
Therefore there exists infinitely many $n$ such that 
$$C_f \rest \bar{I}_n = C_{\vec{f} \rest [v_n,v_{n+1})} \subseteq
T_n.$$
In particular there exists $k \in [v_n,v_{n+1})$ such that 
$\vec{f}(k)=f \rest J_k \in S^k$. Thus there exists $j \in J_k$ such
that 
$$f(j)=s^k(j)=g(j)$$
which finishes the proof of the first part of the theorem.
Note that we only used the the fact that $m(T_n) \leq \frac{1}{4}$ for
$n \in \omega$. 

\vspace{0.1in}

$(2)$ $\unif(\M) \leq \cof({\cal E},\M)$.

To prove this inequality we have to ``dualize'' the above argument. 
Suppose that $ {\cal B} \subseteq \N$ is a family of size ${\lambda}$
witnessing that $\cof({\cal E},\N)={\lambda}$.
We will construct a family $F \subseteq \omega^\omega$ of size
${\lambda}$ such that 
$$\forall f \in \omega^\omega \ \exists g \in F \ \exists^\infty n
\ f(n)=g(n).$$
By \ref{bart}, this will finish the proof.

Since $\cof({\cal E},\N) \geq {\frak b}$ we can find a family $G
\subseteq \omega^\omega$ of size ${\lambda}$ which is unbounded and
consists of increasing functions.

Let $G = \{f_\eta : \eta < {\lambda}\}$ and ${\cal B}=\{H_\eta : \eta
<{\lambda}\}$. 
Without loss of generality we can assume that 
$$H_\eta = \{x \in 2^\omega : \exists^\infty n \ x \rest n \in
H^\eta_n\}$$
where $\sum_{n=1}^\infty |H^\eta_n| \cdot 2^{-n} < \infty$.
For every $\xi,\eta <{\lambda}$ and $n \in \omega$ define
$$\bar{I}^{\xi,\eta}_n = [f_\eta(2n-1),f_\eta(2n+1)]$$
and
$$T^{\xi,\eta}_n = \left\{s \in \bar{I}^{\xi,\eta}_n : \exists j \in
[f_\eta(2n), f_\eta(2n+1)] \ \exists t \in H^\xi_j \ s \rest
\bar{I}^{\xi,\eta}_n = t \rest \bar{I}^{\xi,\eta}_n \right\}$$
Let 
$$W = \left\{\<\xi,\eta\> : \forall n \ |T^{\xi,\eta}_n| \cdot
2^{-|\bar{I}^{\xi,\eta}_n|} \leq 2^{-n} \right\}.$$

Arguing as in the proof of \ref{smalunion}, we show that for every
closed measure zero set $F \subseteq 2^\omega$ there exists
$\<\xi,\eta\> \in W$ such that 
$$\exists^\infty n \ F \rest \bar{I}^{\xi,\eta}_n \subseteq
T^{\xi,\eta}_n .$$

Let $V$ be the set of triples $\<\xi,\eta,\gamma\> \in {\lambda}^3$
such that 
$\<\xi,\eta\> \in W$ and the partition
$\<\lft1[f_\gamma(n),f_\gamma(n+1)\rgt1) : n \in \omega\>$ is finer
that $\<\bar{I}^{\xi,\eta}_n : n \in \omega\>$.

For every triple $\<\xi,\eta,\gamma\> \in V$ let $g^{\xi,\eta,\gamma}
\in \omega^\omega$ be the function $g$ defined in the proof above.

Let 
$$F = \left\{g^{\xi,\eta,\gamma} : \<\xi,\eta,\gamma\> \in
V\right\}.$$
We will show that this family has required properties.
Suppose that $f \in \omega^\omega$.
Find $\gamma,{\delta} < {\lambda}$ such that 
\begin{enumerate}
\item $f_{\delta}(n+1) \geq f_{\delta}(n) +
2^{f_\gamma(n)}\cdot(n+1)$,
\item $\exists^\infty n \ f_\gamma(n+1) > f' \lft1( f_{\delta}(n+1)
\rgt1)^{f_{\delta}(n+1)} + f_\gamma(n)$
\end{enumerate}
where $f'(n)=\max\{f(1), \ldots, f(n)\}+1$.

Define
$I_n= \lft1[f_\gamma(n),f_\gamma(n+1) \rgt1) $ and $J_n=
\lft1[f_{\delta}(n),f_{\delta}(n+1) \rgt1) $ for $n \in \omega$.
As in the above part we have
$$\exists^\infty n \ \vec{f}(n) \in I_n.$$
Now we can find $\<\xi,\eta\> \in W$ such that 
$$\exists^\infty n \ C_f \rest \bar{I}^{\xi,\eta}_n \subseteq
T^{\xi,\eta}_n .$$
It follows that
$$\exists^\infty n \ f(n)=g^{\xi,\eta,\gamma}(n)$$
which finishes the proof.~$\QED$

We conclude this section with two applications.

In [Mi1] it is proved that:
\begin{theorem}[Miller]
$\add(\N) \leq {\frak b}$ and $\cof(\N) \geq {\frak d}$.~$\QED$
\end{theorem}

\begin{theorem}[Bartoszynski, Raisonnier, Stern {[Ba], [RS]}]
$\add(\N) \leq \add(\M)$ and $\cof(\N) \geq \cof(\M)$.
\end{theorem}
\Proof
We have
$$\add(\N) \leq \min\{{\frak b}, \add({\cal E},\N)\} = \min\{{\frak b},
\cov(\M)\}= \add(\M).$$
Similarly
$$\cof(\N) \geq \max\{{\frak d},\cof({\cal E},\N)\}=\max\{{\frak d},
\unif(\M)\} = \cof(\M).~\QED$$

Also we get another proof of the main result from [BJ]:
\begin{theorem}[Bartoszynski, Judah]
$\cf \lft1( \cov(\M) \rgt1) \geq \add(\N)$.
\end{theorem}
\Proof
Clearly $\cf \lft1( \add({\cal E},\N) \rgt1) \geq \add(\N)$.~$\QED$

\section{Cardinals $\cov({\cal E})$ and $\unif({\cal E})$}
\setcounter{theorem}{0}

In this section we will prove some results concerning covering number
of ${\cal E}$. Most of the results are implicite in  [Ba2] and [BJ1].

Let us start with the following easy observation.

\begin{lemma}\label{1520}
\ \begin{enumerate}
\item Every null set can be covered by ${\frak d}$ many closed null sets,
\item Every null set of size $< {\frak b}$ can be covered by a null set of type
$F_{\sigma}$.
\end{enumerate}
\end{lemma}
\Proof
Suppose that $G$ is a null subset of $2^{\omega}$. As in \ref{oxtoby},
we can assume that
$$G = \left\{x \in 2^{\omega} : \exists^{\infty} n \ x \rest n \in F_{n}\right\}$$
where 
$\sum_{n=1}^{\infty} |F_{n}| \cdot 2^{-n} < \infty$.
For every $x \in G$ let $f_{x} \in \omega^{\omega}$  be an increasing
enumeration of the set $\{n \in \omega: x \rest n \in F_{n}\}$.
For a strictly increasing function $f \in \omega^{\omega}$ let
$$G_{f} = \left\{x \in 2^{\omega}: \forall^{\infty}n
\ \exists m \in [n, f(n)] \ x \rest m \in F_{m}\right\}  .$$
It is clear that for every
$f  \in \omega^{\omega}$ the set $G_{f} \subseteq G$ is a measure zero  set
of type $F_{\sigma}$.
Notice also that if $f_{x} \leq^{\star} f$ then $x \in G_{f}$.

\vspace{0.1in}

$(1)$ Let $F \subseteq \omega^{\omega}$ be a dominating family of size
${\frak d}$ which consists of increasing functions.
Then by the above remarks
$$G = \bigcup_{f \in F} G_{f}  .$$

$(2)$ Suppose that $X \subseteq G$ is a set of size $< {\frak b}$.
Let  $f$ be an increasing function which dominates all functions
$\left\{f_{x} : x \in X\right\}$. Then $X \subseteq G_{f}$.~$\QED$

As a corollary we get:
\begin{theorem}
\ \begin{enumerate}
\item  If $\cov({\cal M}) = {\frak d}$ then
$\cov({\cal E}) = \max\left\{\cov({\cal M}), \cov({\cal N})\right\}$,
\item  If $\unif({\cal N}) = {\frak b}$ then
$\unif({\cal E}) = \min\left\{\unif({\cal M}), \unif({\cal N})\right\}$.
\end{enumerate}
\end{theorem}
\Proof
Since  ${\cal E} \subseteq {\cal M} \cap {\cal N}$
we have
$$\cov({\cal E}) \geq
\max\left\{\cov({\cal M}), \cov({\cal N})\right\}$$
and
$$\unif({\cal E}) \leq \min\left\{\unif({\cal M}), \unif({\cal N})\right\} .$$
By the previous lemma
$$\max\left\{ \cov(\M), {\frak d}\right\} \geq \cov({\cal E})$$ and
$$\unif({\cal E}) \geq \min\left\{ \unif(\N), {\frak b}\right\} $$
which finishes the proof.~$\QED$

Suppose that $f \in \omega^{\omega}$ and
$\sum_{n=1}^{\infty} 2^{-f(n)} < \infty$.
Define
$$\varphi \in \Sigma_{f}
\leftrightarrow \varphi \in \lft1( [\omega]^{<\omega}\rgt1)^{\omega} \&
\ \forall n \ \left(\varphi(n) \subseteq 2^{f(n)}  \
\& \ 
\frac{|\varphi(n)|}{2^{f(n)}} \leq \frac{1}{4^n}\right) $$
and
$$\varphi  \in \Pi_{f}
\leftrightarrow \varphi  \in \lft1( [\omega]^{<\omega}\rgt1)^{\omega}\ \&
\ \forall n \ \varphi(n) \subseteq 2^{f(n)} \ 
\& \ \exists^\infty n \ 
\frac{|\varphi(n)|}{2^{f(n)}} \leq \frac{1}{4^n} $$ 
and let
$ {\cal X}_{f} = \prod_{n=1}^{\infty} 2^{f(n)}  .$

Notice that  $\Sigma_{f} \subseteq \Pi_{f}$.

For $\varphi \in \Sigma_{f} \cup \Pi_f$ define
define set $H_{\varphi} \subseteq 2^{\omega}$
as follows:

Let $k_{n} = 1+2+ \cdots + f(n)$ for $n \in \omega$.
Identify natural numbers 
$\leq 2^{f(n)}$ with 0-1 sequences of length $f(n)$ and
define
$$H_{\varphi} = \left\{x \in 2^{\omega} :\  \forall^{\infty}n \ x \rest
[k_{n},
k_{n+1}) \in
\varphi(n) \right\}  .$$
Note that
$$\mu(H_{\varphi}) \leq \prod_{n=m}^{\infty} \mu(\left\{x \in 2^{\omega}:
x \rest [k_{n}, k_{n+1}) \in \varphi(n)\right\}) \leq
\sum_{m=1}^{\infty} \prod_{n=m}^{\infty} \frac{|\varphi(n)|}{2^{f(n)}}
= 0 .$$

For $x \in 2^{\omega}$. Define $h_{x}(n) = x \rest [k_{n}, k_{n+1})$
for $n \in \omega$.
Clearly $h_{x}$ corresponds to an element of ${\cal X}_{f}$.

Finally we have
$$x \in H_{\varphi}\  \leftrightarrow  \ \forall^{\infty} n \ h_{x}(n) \in
\varphi(n) . $$

\begin{theorem}\label{158}
Suppose  that $C \in {\cal E}$. Then  there exists $f \in
\omega^\omega$ and $\varphi \in \Sigma_f$ such that 
$C \subseteq H_\varphi $.
\end{theorem}
\Proof
Suppose that $C \subseteq 2^{\omega}$ is a null set of type
$F_{\sigma}$. Represent $C$ as $\bigcup_{n \in \omega} C_{n}$ where
$\<C_{n} : n \in \omega\>$ is an increasing family of closed
sets of measure zero.
Define sequence $\<k_{n} : n \in \omega\>$ as follows: $k_{0} = 0$ and
$$k_{n+1} = \min \left\{m > k_{n}: \exists \ T_{n} \subseteq 2^{m} \ \left(
C_{n} \subseteq [T_n] \ \& \ \frac{|T_{n}|}{2^{m}}
\leq \frac{1}{4^{k_{n}}}\right)\right\} .$$
Let $I_{n} = [k_{n}, k_{n+1})$ and $J_{n} = \left\{s \rest I_{n} : s \in T_{n} \right\}$
for $n \in \omega$.
We can see that for all $n \in \omega$
$$\frac{|J_{n}|}{2^{|I_{n}|}} \leq 2^{k_{n}} \cdot \frac{1}{4^{k_{n}}} \leq
\frac{1}{2^{n}}  .$$
We also have
$$F \subseteq \left\{x \in 2^{\omega} :  \forall^{\infty} n \ x \rest I_{n} \in
J_{n} \right\} =H_\varphi$$
where
$f(n)=|I_n|$ and $\varphi(n)=J_n$ for all $n$. By the above remarks
$\varphi \in \Sigma_f$.~$\QED$

For an increasing function $g \in \omega^\omega$ define $g^\star \in
\omega^\omega$ as $g^\star(0)=0$ and $g^\star(n+1)=g(g^\star(n)+1)$.
\begin{lemma}\label{simp}
Suppose that $f,g \in \omega^\omega$ are increasing functions and
$\varphi \in \Sigma_f$.
\begin{enumerate}
\item If $f \leq^\star g$ then there exists $\psi \in
\Sigma_{g^\star}$ such that $H_\varphi \subseteq H_\psi$,
\item if $g \not \leq^\star f$  then there exists $\psi \in
\Pi_{g^\star}$ such that $H_\varphi \subseteq H_\psi$.
\end{enumerate}
\end{lemma}
\Proof
Let $I_n=\lft1[f(n),f(n+1)\rgt1) $ and $I^\star_n =
\lft1[g^\star(n),g^\star(n+1)\rgt1) $
for $n \in \omega$.
Note that if $f \leq^\star g$ then
$$\forall^\infty n \ \exists m \ I_m \subseteq I^\star_n$$
and if $g \not \leq^\star f$ then
$$\exists^\infty n \ \exists m \ I_m \subseteq I^\star_n .$$
Define 
$$\psi(n) = \left\{ \begin{array}{ll}
\left\{s \in 2^{I^\star_n} : 
\exists m \ \lft1( I_m \subseteq I^\star_n \ \& \ s \rest I_m \in
\varphi(m)\rgt1) \right\}
& \hbox{if } \exists m \ I_m  \subseteq I^\star_n\\
2^{I^\star_n} & \hbox{otherwise}
\end{array}
\right. $$
It follows that $\psi \in \Sigma_{g^\star}$ in the first case and
$\psi \in \Pi_{g^\star}$ in the second case.
Moreover, the inclusion,
$H_\varphi \subseteq H_\psi$ is an immediate consequence of
the above definition.~$\QED$

As a consequence we get:
\begin{theorem}\label{159a}
Suppose that $\left\{F_\xi : \xi < \kappa\right\}$ is a family of elements of
${\cal E}$. 
\begin{enumerate}
\item If $\kappa  < {\frak b}$ then there exists a function $g \in
\omega^\omega$ and a family $\left\{\varphi_\xi :\xi < \kappa \right\} \subseteq
\Sigma_g$ such that $F_\xi \subseteq H_{\varphi_\xi}$ for $\xi <
\kappa$,
\item if $\kappa  < {\frak d}$ then there exists a function $g \in
\omega^\omega$ and a family $\left\{\varphi_\xi :\xi < \kappa \right\} \subseteq
\Pi_g$ such that $F_\xi \subseteq H_{\varphi_\xi}$ for $\xi < \kappa$. $\QED$
\end{enumerate}
\end{theorem}

The following fact follows immediately from \ref{159a}.
\begin{theorem}\label{159b}
If\/ $\cov({\cal E}) < {\frak d}$ then there exists $f \in \omega^{\omega}$
such that $\cov({\cal E})$ is equal to the size of the smallest family
$\Psi \subseteq \Pi_{f}$ such that
$$\forall h \in {\cal X}_{f} \ \exists \psi \in \Psi \ \forall^{\infty} n
\ h(n) \in \psi(n)  . \ \QED$$
\end{theorem}

As an corollary we get the following:
\begin{theorem}[Miller]
If $\cov({\cal E}) \leq {\frak d}$ then $\cf \lft1(\cov({\cal E})\rgt1) >\aleph_{0}$.
\end{theorem}
\Proof
Suppose that $\cf(\cov({\cal E})) = \aleph_{0}$. 
Since ${\frak d}$ has uncountable
cardinality we have
$\cov({\cal E}) < {\frak d}$.
By \ref{159b}
under this assumptions there exists $g \in \omega^{\omega}$
such that
$\cov({\cal E})$ is the size of the smallest family
$\Psi \subseteq \Pi_{g}$ such that
$$\forall h \in {\cal X}_{g} \ \exists \psi \in \Psi \ \forall^{\infty}
n \ h(n) \in \psi(n)  .$$
Assume that $\Psi$ is the smallest family having above properties and
let $\left\{\Psi_{n} : n \in \omega\right\}$ be an increasing family such that
$\Psi = \bigcup_{n \in \omega} \Psi_{n}$ and $|\Psi_{n}| < |\Psi|$ for
all $n \in \omega$.

By the assumption for every $m \in \omega$ there exists a function
$h_{m} \in {\cal X}_{g}$ such that
$$\forall m \ \forall \psi \in \Psi_{m} \ \exists^{\infty} n
\ h_{m}(n) \not \in \psi(n)  .$$
For $\psi \in \Psi$ define  $k^{\psi}_{0} = 0$ and for $n \in \omega$
$$k^{\psi}_{n+1} = \min\left\{m > k^{\psi}_{n} : \forall j \leq n \
\exists i \in [k^{\psi}_{n},m) \ h_{j}(i) \not \in \psi(i)\right\} .$$
Since $|\Psi| < {\frak d}$ we can find an increasing function
$r \in \omega^{\omega}$ such that 
$$\forall \psi \in \Psi  \ \exists^{\infty} n \
k^{\psi}_{n} \leq r(n)  .$$
Let $h = h_{1} \rest [r^\star(0),r^\star(1))^{\frown}h_{2} 
\rest [r^\star(1),r^\star(2))^{\frown}
h_{3} \rest [r^\star(2),r^\star(3))^{\frown} \ \ldots$ .

Fix $\psi \in \Psi$. By the assumption about $r$ we have
$$\exists^{\infty}n \ \exists m > n \ 
r^\star(n) < k^{\psi}_{m} < k^{\psi}_{m+1} <
r^\star(n+1) .$$
But this means that
$$\exists i \in [r(n),r(n+1)) \ h_{n+1}(i) = h(i) \not \in \psi(i) .$$
Since $\psi$ is an arbitrary element of $\Psi$ it finishes the proof.~$\QED$

\section{Consistency results}
\setcounter{theorem}{0}
The goal of this section is to show that $\cov({\cal E}) > \max\left\{
\cov(\N), \cov(\M)\right\}$ and 
$\unif({\cal E})< \min\left\{\unif(\N), \unif(\M)\right\}$ are both consistent
with ZFC.  We use the technique developed in [JS].

\begin{lemma}\label{umph}
Suppose that ${\cal P}$ is a notion of forcing satisfying ccc.
Let $\dot{C}$ be a ${\cal P}$-name for an element of ${\cal E}$. 
\begin{enumerate}
\item If ${\cal P}$ does not add dominating reals
then there exists $f \in \omega^\omega \cap \V$ and a ${\cal P}$-name
$\dot{\varphi}$ such that\/ $\forces_{{\cal P}} \dot{\varphi} \in \Pi_f$
and\/ $\forces_{{\cal P}} \dot{C} \subseteq H_{\dot{\varphi}}$,
\item  if ${\cal P}$ is $\omega^\omega$-bounding
then there exists $f \in \omega^\omega \cap \V$ and a ${\cal P}$-name
$\dot{\varphi}$ such that\/ $\forces_{{\cal P}} \dot{\varphi} \in \Sigma_f$
and\/ $\forces_{{\cal P}} \dot{C} \subseteq H_{\dot{\varphi}}$.
\end{enumerate}
\end{lemma}
\Proof Follows immediately from \ref{simp}. $\QED$

\begin{definition}\label{big}
Suppose that $N \models ${\rm ZFC}$^\star$. A function $x \in 2^\omega$ is
called $N$-big iff 
$$x \not \in \bigcup ({\cal E} \cap N) .$$

We say that a partial ordering ${\cal P}$ satisfying ccc is good
if for every model 
$N \prec H(\chi)$ and every filter $G$ which is ${\cal
P}$-generic over $\V$,  if $x \in 2^\omega$ is $N$-big then $x$ is $N[G]$-big.
\end{definition}

Let ${\bf B}$ denote the random real forcing.
\begin{theorem}
${\bf B}$ is good.
\end{theorem}
\Proof
Suppose that $x$ is $N$-big. Let $\dot{C} \in N$ be a ${\bf B}$-name for an
element of~${\cal E}$. Since ${\bf B}$ is $\omega^\omega$-bounding, by
\ref{umph}, we can find a function 
$f  \in \omega^\omega \cap N$ and a ${\bf B}$-name $\dot{\varphi} \in
N$ for an element of $\Sigma_f$ such that $\forces_{{\newbf B}}
\dot{C} \subseteq H_{\dot{\varphi}}$.

For $s \in 2^{f(n)}$ define $B_{n,s} = \lbv s \in \dot{\varphi}(n)
\rbv_{{\newbf B}}$. 
Let 
$$\varphi(n) = \left\{s : \mu(B_{n,s}) \geq \frac{1}{2^n}\right\} \hbox{ for } n \in
\omega . $$
Note that since 
$$\forces_{{\newbf B}} \frac{|\dot{\varphi}(n)|}{2^{f(n)}} \leq \frac{1}{4^n}$$
we get that
$$\frac{|\varphi(n)|}{2^{f(n)}} \leq \frac{1}{2^n} \hbox{ for } n \in \omega.$$

Suppose that $p \forces_{{\newbf B}} \forall n \geq m \ x \rest n \in
\dot{\varphi}(n)$. Find $k$ such that $\mu(p) \geq 2^{-k}$. Since $x$ is
$N$-big there exists $n \geq k$ such that $\widehat{s}=x \rest I_n \not \in
\varphi(n)$. In particular $\mu(B_{n,\hat{s}}) < 2^{-k}$. 
Let $q = p - B_{n,\hat{s}}$.
It is clear that
$$q \forces_{{\newbf B}} x \rest I_n \not \in \dot{\varphi}(n)$$
which gives a contradiction.~$\QED$

\begin{lemma}\label{pres}
\
\begin{enumerate}
\item If ${\cal P}$ and ${\cal Q}$ are good forcing notions then ${\cal P}
\star {\cal Q}$ is good.
\item If $\left\{{\cal P}_\alpha , \dot{{\cal Q}}_\alpha : \alpha < \delta\right\}$ is
a finite support iteration such that
\begin{enumerate}
\item $\forces_\alpha \dot{{\cal
Q}}_\alpha \hbox{ is good} $,
\item $\forces_\alpha \hbox{``} \omega^\omega \cap \V \hbox{ is
unbounded''}.$
\end{enumerate}
then ${\cal P}_\delta = \lim_{\alpha <
\delta} {\cal P}_\alpha $ is good.
\end{enumerate}
\end{lemma}
\Proof
The first part is obvious.
We will prove the second part by induction on $\delta$. 
Without loss of generality we can assume that $\delta$ is a limit
ordinal. Suppose that the lemma is true for $\alpha<\delta$. Let $N
\prec H(\chi)$
be a model and let 
$\dot{C}$ be a ${\cal P}_\delta$-name for an element of ${\cal E} \cap N$.
It is well known that under the assumptions  
${\cal P}_\delta$ does not add dominating reals. Therefore there exists $f
\in \omega^\omega \cap N$ and a ${\cal P}_\delta$-name
$\dot{\varphi}$ for an element of $\Pi_f$ such that 
$$\forces_\delta \dot{C} \subseteq H_{\dot{\varphi}} .$$
Assume that $x$ is $N$-big and suppose that for some $p \in {\cal
P}_\delta$,
$$p \forces_\delta \forall n > n_0 \ x \rest [f(n),f(n+1)) \in
\dot{\varphi}(n). $$
Define a sequence $\<p_n : n \in \omega\>, \<k_n: n \in \omega\> \in N$ and
$\varphi \in \Pi_f$ such that 
\begin{enumerate}
\item $p = p_0 \leq p_1 \leq p_2 \ldots $,
\item $p_{n+1} \forces_\delta \forall j \leq k_n \ \dot{\varphi}(j) =
\varphi(j)$,
\item $p_{n+1} \forces_{\delta} \exists j \in [k_n,k_{n+1}] \
|\varphi(j)| \cdot 2^{f(j)} \leq 4^{-j}$.
\end{enumerate}
Since $x$ is $N$-big there exists $m > n_0$ such that $x \rest
[f(m),f(m+1)) \not \in \varphi(m)$. Therefore
$p_{m} \forces x \rest [f(m),f(m+1)) \not \in
\dot{\varphi}(m)$. In particular,
$$p \not \forces x \rest [f(m),f(m+1)) \not \in
\dot{\varphi}(m)$$
which is a contradiction.~$\QED$

\begin{theorem}
It is consistent with {\rm ZFC} that 
$$\unif({\cal E}) <
\min\left\{\unif(\N), \unif(\M)\right\}.$$
\end{theorem}
\Proof
Let ${\cal P}_{\omega_2}$ be a finite support iteration of length
$\omega_2$ of random real forcing. Let $G$ be a ${\cal
P}_{\omega_2}$-generic filter over a model $\V \models GCH$.
Since ${\cal P}_{\omega_2}$ adds random and Cohen reals we have
$\V[G] \models \unif(\M)=\unif(\N) = \aleph_2$.
We will show that $\V[G] \models \unif({\cal E})=\aleph_1$. It is
enough to show that $\V[G] \models 2^\omega \cap \V \not \in {\cal
E}$.

Suppose that $C \in \V[G] \cap {\cal E}$. Let $\dot{C}$ be a ${\cal
P}_{\omega_2}$-name for $C$. Let $N \prec H(\chi)$ be a countable
model containing $\dot{C}$ and ${\cal P}_{\omega_2}$.
Since $N$ is countable there exists $x \in 2^\omega \cap \V$ which is
$N$-big. By \ref{pres}, $x$ is also $N[G]$-big. In particular $x
\not \in C$.~$\QED$

\begin{theorem}
It is consistent with {\rm ZFC} that $\cov({\cal E}) >
\max\left\{\cov(\N), \cov(\M)\right\}$.
\end{theorem}
\Proof
Let ${\cal P}_{\omega_1}$ be a finite support iteration of length
$\omega_1$ of random real forcing. Let $G$ be a ${\cal
P}_{\omega_1}$-generic filter over a model $\V \models \cov({\cal
E})=\aleph_2$.

It is clear that $\V[G] \models \cov(\N) = \cov(\M)=\aleph_1$. We will
show that $\V[G] \models \cov({\cal E})=\aleph_2$.

Suppose that $\left\{C_\xi : \xi < \omega_1\right\} \subseteq \V[G] \cap {\cal
E}$. Let $\dot{C}_\alpha$ be a ${\cal P}_{\omega_1}$-name for
$C_\alpha$.
Let $N \prec H(\chi)$ be a model of size $\aleph_1$ containing all
names $\dot{C}_\alpha$ and ${\cal P}_{\omega_1}$. Since $\V \models
\cov({\cal E}) > \aleph_1$ there exists $x \in 2^\omega \cap \V$ which
is $N$-big. By \ref{pres}, $x$ is also $N[G]$-big. In particular, $x
\not \in \bigcup_{\xi < \omega_1} C_\xi$.~$\QED$

\end{document}